\documentclass{article}
\usepackage{amsmath,amssymb,amsfonts}

\allowdisplaybreaks[3]

\newcommand{\ncd}{\newcommand}

\def\R{\mathbb R}
\def\Z{\mathbb Z}
\def\BEQ#1\EEQ{\begin{equation}#1\end{equation}}
\def\BAL#1\EAL{\begin{align}#1\end{align}}
\ncd{\gvct}[1]{{\boldsymbol{#1}}}
\ncd{\bma}{\begin{pmatrix}}
\ncd{\ema}{\end{pmatrix}}
\ncd{\privateremark}[1]{}
\newtheorem{theo}{Theorem}[section]
\newtheorem{dfn}{Definition}[section]
\newtheorem{lemma}[theo]{Lemma}

\DeclareMathOperator{\glcd}{glcd}

\begin{document}

\title{Multiple CSLs for the body centered cubic lattice}

\author{Peter Zeiner\\
Faculty of Mathematics, Bielefeld University\\ 
Universit\"atsstra{\ss}e 25, 33615 Bielefeld, Germany\\
and\\
Institute of Theoretical Physics \& CMS, TU Wien\\
Wiedner Hauptstra{\ss}e 8--10, 1040 Wien, Austria}

\maketitle

\begin{abstract}
Ordinary Coincidence Site Lattices (CSLs) are defined as the intersection
of a lattice $\Gamma$ with a rotated copy $R\Gamma$ of itself. They are useful
for classifying grain boundaries and have been studied extensively since the
mid sixties. Recently the interests turned to so-called multiple CSLs, i.e.
intersections of $n$ rotated copies of a given lattice $\Gamma$, in particular
in connection with lattice quantizers. Here we consider multiple CSLs
for the 3-dimensional body centered cubic lattice. We discuss
the spectrum of coincidence indices and their multiplicity, in particular we
show that the latter is a multiplicative function and give an explicit
expression of it for some special cases.
\end{abstract}

\section{Introduction}

Ordinary coincidence site lattices (CSLs) have been studied intensively 
since the 1960s
(see e.g. \cite{boll70,boll82} and references therein), because they are an
important tool to characterize and
analyze the structure of grain boundaries in crystals. Hence there is a vast
literature on 3--dimensional CSLs, in particular on cubic
lattices~\cite{grim74,gribo74,ble81,grim84,baa97,pzcsl1}. CSLs in higher
dimensions have been studied as well (see~\cite{baa97,pzcsl2} and references
therein) with possible applications to quasiperiodic structures. In addition,
the concept of CSLs has been generalized for modules, again to cover the needs
of quasiperiodic structures~\cite{plea96,baa97}.

We want to discuss another generalization of CSLs here. Loosely speaking, CSLs
are the intersection of two mutually rotated lattices. It is thus natural to
consider the intersection of $n$ mutually rotated lattices. This question has
recently been raised in connection with quantizing
procedures~\cite{slo02a,slo02b}, where it seems useful to represent a complex
lattice as the intersection of simpler lattices.
In the meantime this question has been answered in detail for several
2--dimensional lattices~\cite{baagri05}, so the next step is to consider a
$3$--dimensional example. We choose a cubic lattice here, namely the
body--centered cubic lattice,
since the cubic case is one of the best studied $3$--dimensional cases. It
might seem more natural to discuss the primitive cubic case first, but it
turns out that the body centered case can be treated more elegantly. Moreover
most results hold for all three cubic lattices. Details for the other cubic
lattices shall be published elsewhere~\cite{pzmcsl2}.

\section{Ordinary CSLs}

We recall some definitions first~\cite{baa97}. Let
$\Gamma\subseteq\R^n$ be
an $n$-dimensional lattice and $R\in SO(n)$ a rotation (we restrict our
considerations here to proper rotations for simplicity). Then $R$ is called a
\emph{coincidence rotation} if $\Gamma(R):=\Gamma\cap R\Gamma$ is a lattice of
finite index in $\Gamma$. The corresponding lattice
$\Gamma(R)=\Gamma\cap R\Gamma$ is called a \emph{coincidence site lattice}.
The \emph{coincidence index} $\Sigma(R)$ is defined as the index of
$\Gamma(R)$ in $\Gamma$. Note in passing that the set of coincidence rotations
forms a group under matrix multiplication.

In the following let $\Gamma$ be a body centered cubic lattice.
Then $R\in SO(3)$ is
a coincidence rotation if and only if $R$ is a matrix with rational entries
(see~\cite{grim74,gribo74,ble81,grim84,baa97,pzcsl1}). They can be
parameterized by integral quaternions $\gvct{q}=(\kappa,\lambda,\mu,\nu)$,
i.e., by quaternions with integral coefficients $\kappa,\lambda,\mu,\nu$ in the
following way~\cite{koecheng,hurw,val,baa97,pzcsl1}:
\BEQ
R(\gvct{q})=\frac{1}{|\gvct{q}|^2}
\begin{pmatrix}
\kappa^2+\lambda^2-\mu^2-\nu^2 & -2\kappa\nu+2\lambda\mu &
2\kappa\mu+2\lambda\nu\\
2\kappa\nu+2\lambda\mu & \kappa^2-\lambda^2+\mu^2-\nu^2 &
-2\kappa\lambda+2\mu\nu\\
-2\kappa\mu+2\lambda\nu & 2\kappa\lambda+2\mu\nu &
\kappa^2-\lambda^2-\mu^2+\nu^2 
\end{pmatrix},
\EEQ
where $|\gvct{q}|^2=\kappa^2+\lambda^2+\mu^2+\nu^2$ is called the norm of
$\gvct{q}$. Note that we will call a quaternion an integer
quaternion if it is an integral quaternion $(\kappa,\lambda,\mu,\nu)$
or the sum of an integral quaternion with the quaternion $(1,1,1,1)/2$.
We call an integral quaternion $\gvct{r}=(\kappa,\lambda,\mu,\nu)$ primitive
if the greatest common divisor of $\kappa,\lambda,\mu,\nu$ equals $1$. If
not stated otherwise every (integral) quaternion will be assumed to be a
primitive quaternion. Furthermore let $\gvct{H}\subset\gvct{Q}$ denote the
ring of integer quaternions and the ring of real quaternions, respectively.

One can show that the coincidence index is given by
$\Sigma(R(\gvct{q}))=|\gvct{q}|^2/2^\ell$, where $\ell$ is the maximal power
such that $2^\ell$ divides $|\gvct{q}|^2$
(see e.g.~\cite{gribo74,grim84,baa97}), i.e. $\Sigma(R)$ is always odd. On the
other hand $\Sigma(R)$ runs over all positive odd integers if $R$ runs over
all coincidence rotations, i.e. the spectrum of coincidence rotations is the
set of all positive odd integers. Let $O$ denote the cubic symmetry
group. Then $\Gamma(RQ)=\Gamma(R)$ if and only if $Q\in O$, i.e.
$R$ and $RQ$ generate the same CSL, which motivates to call $R$ and $RQ$
\emph{strongly equivalent}. In general one calls $R$ and $R'$
\emph{equivalent} if there exist $Q,Q'$ such that $R=QR'Q'$. For any
$R(\gvct{q})$ we can find a strongly equivalent $R'=R(\gvct{q}')$ such that
$|\gvct{q}'|^2$ is odd, i.e. $\Sigma(R')=|\gvct{q}'|^2$.
Hence we will assume in the following that $|\gvct{q}|^2$ is odd.

If we define the projection $P:\gvct{Q}\to \R^3$ by
$P(q_0,q_1,q_2,q_3)=(q_1,q_2,q_3)$ then $\Gamma=P\gvct{H}$, i.e. the body
centered cubic lattice is obtained by a projection of $\gvct{H}$ onto $\R^3$.
Moreover Lemma 5.1 of~\cite{pzcsl1} tells us that
$\Gamma(R(\gvct{q}))=P(\gvct{q}\gvct{H})$, where $\gvct{q}\gvct{H}$ is a left
ideal of $\gvct{H}$. In fact this establishes a one to one correspondence of
CSLs and left ideals of $\gvct{H}$ (see~\cite{pzmcsl2}), which is a key
in the discussion of the body centered CSLs. In particular, finding the number
$f(\Sigma)$ of different CSLs of given index $\Sigma$ is equivalent to
counting the corresponding left ideals of $\gvct{H}$. One can show that
$f(\Sigma)$ is a multiplicative function, i.e. $f(mn)=f(m)f(n)$ if $m$ and $n$
are coprime, and in particular we have $f(1)=1$, $f(2)=0$
and $f(p^r)=(p+1)p^{r-1}$ for all odd primes $p$~\cite{baa97,pzcsl1}.

\section{Multiple CSLs}

Now we turn to multiple CSLs, which we define as
follows~\cite{baagri05,pzmcsl2}:
\begin{dfn}
Let $\Gamma$ be an $n$-dimensional lattice and $R_i$, $i=1,\ldots m$
coincidence rotations of $\Gamma$. Then the lattice
\BAL
\Gamma(R_1,\ldots,R_m):=\Gamma\cap R_1\Gamma\cap\ldots\cap R_m\Gamma=
\Gamma(R_1)\cap\ldots\cap\Gamma(R_m)
\EAL
is called a multiple CSL (MCSL). Its index in $\Gamma$ is denoted by
$\Sigma(R_1,\ldots,R_m)$.
\end{dfn}
Note that $\Sigma(R_1,\ldots,R_m)$ is finite since $\Gamma(R_1,\ldots,R_m)$
is a finite intersection of mutually commensurate lattices~\cite{baa97}. In
particular, it follows from the second homomorphism theorem that
\BAL
\Sigma(R_1,R_2)=\frac{\Sigma(R_1)\Sigma(R_2)}{\Sigma_+(R_1,R_2)},
\EAL
where $\Sigma_+(R_1,R_2)$ is the index of the direct sum
$\Gamma_+(R_1,R_2)=\Gamma(R_1)+\Gamma(R_2)$ in $\Gamma$. In general one shows
\BAL
\Sigma(R_1,\ldots,R_m)=\frac{\Sigma(R_1,\ldots,R_{m-1})\Sigma(R_m)}%
{\Sigma_+(R_1,\ldots,R_{m-1};R_m)},
\EAL
where $\Sigma_+(R_1,\ldots,R_{m-1};R_m)$ is the index of
$\Gamma_+(R_1,\ldots,R_{m-1};R_m)=\Gamma(R_1,\ldots,R_{m-1})+\Gamma(R_m)$
in $\Gamma$. In particular, $\Sigma(R_1,\ldots,R_m)$ divides
$\Sigma(R_1)\cdot\ldots\cdot\Sigma(R_m)$. In case of the 3--dimensional
cubic lattices this implies immediately that the spectrum is the same
for multiple and ordinary CSLs, i.e. $\Sigma(R_1,\ldots,R_m)$ runs over all
odd positive integers, too.
However, new lattices emerge and the
multiplicity of a given index will increase. Note that the spectrum is
preserved also for the square lattice and the 4--dimensional hypercubic
lattices. 

Having determined the spectrum we can attack the second main problem, the
number $f_m(\Sigma)$ of different MCSLs. To this end we have to
determine all possible MCSLs. We first note that
$\Gamma_+(R_1,R_2)=P(\gvct{q}_1\gvct{H}+\gvct{q}_2\gvct{H})=
P(\gvct{q}\gvct{H})$, where $\gvct{q}$ is the greatest left common divisor 
(glcd) of $\gvct{q}_1$ and $\gvct{q}_2$. Hence (recall that we may assume
that $|\gvct{q}_i|^2$ is odd)
\BAL
\Sigma(R_1,R_2)&=\frac{|\gvct{q}_1|^2|\gvct{q}_2|^2}{|\gvct{q}|^2}
\quad \mbox{with }\gvct{q}=\glcd(\gvct{q}_1,\gvct{q}_2).
\EAL
In case that $|\gvct{q}_1|^2$ and $|\gvct{q}_2|^2$ are relatively prime
this reduces to $\Sigma(R_1,R_2)=|\gvct{q}_1|^2|\gvct{q}_2|^2$ which suggests
the following lemma~\cite{pzmcsl2}:
\begin{lemma}
Let $|\gvct{q}_1|^2$ and $|\gvct{q}_2|^2$ be relatively prime. Then there
exists a quaternion $\gvct{q}$, such that $\Gamma(R_1,R_2)=
\Gamma(R(\gvct{q}))$.
\end{lemma}
The proof makes use of the fact there exists a right least common multiple
$\gvct{q}$ of $\gvct{q}_1$ and $\gvct{q}_2$. Note that in this case the MCSL
is equal to an ordinary CSL. This result can be immediately generalized for
arbitrary $m$. Conversely we have~\cite{pzmcsl2}
\begin{lemma}
Let $p_1^{\alpha_1}\cdots p_\ell^{\alpha_\ell}$ be the prime decomposition of
$|\gvct{q}|^2$. Then there exist quaternions $\gvct{q}_i$ such that
$|\gvct{q}_i|^2=p_i^{\alpha_i}$ and
$\Gamma(R(\gvct{q}))=
\Gamma(R(\gvct{q}_1))\cap\ldots\cap\Gamma(R(\gvct{q}_\ell))$.
\end{lemma}
This decomposition is unique. More generally we have~\cite{pzmcsl2}
\begin{lemma}
For any $\Gamma(R_1,\ldots,R_m)$ with
$\Sigma(R_1,\ldots,R_m)=p_1^{\alpha_1}\cdots p_\ell^{\alpha_\ell}$
(all $p_i$ distinct)
there exists a unique decomposition
$\Gamma(R_1,\ldots,R_m)=\Gamma_1\cap\ldots\cap\Gamma_\ell$ such that
$\Gamma_i$ is an MCSL with index $\Sigma_i=p_i^{\alpha_i}$.
\end{lemma}
Thus the analysis of MCSLs can be reduced to the study of the MCSLs with prime
power index. Another consequence is the multiplicativity of $f_m$:
\begin{theo}
Let $f_m(\Sigma)$ be the number of different MCSLs for a given index $\Sigma$.
Then $f_m(nn')=f_m(n)f_m(n')$ if $n$ and $n'$ are relatively prime. 
\end{theo}
Since the analysis of general MCSLs with prime power index is rather
cumbersome we confine our discussion to the intersection of two CSLs. The
general discussion will be published elsewhere~\cite{pzmcsl2}.  We may confine
our discussion 
to the case that neither $\gvct{q}_1$ nor $\gvct{q}_2$ is a right multiple
of the other, i.e. there does not exist an integer quaternion $\gvct{r}$
such that $\gvct{q}_1=\gvct{q}_2\gvct{r}$ or vice versa, since otherwise
$\Gamma(R_1,R_2)$ reduces to $\Gamma(R_2)$ or $\Gamma(R_1)$, respectively.
We first mention a representation of $\Gamma(R_1,R_2)$.
Here
$\gvct{\bar{q}}=(\kappa,-\lambda,-\mu,-\nu)$ denotes
the conjugate of $\gvct{q}=(\kappa,\lambda,\mu,\nu)$.
\begin{lemma}
Let $|\gvct{q}_i|^2=p^{\alpha_i}$, $i=1,2$, $p$ prime, none of the 
$\gvct{q}_i$ a right
multiple of the other one. Choose $\gvct{r}$
such that $\gvct{q}_1\gvct{r}\gvct{\bar{q}}_2$ is a primitive quaternion
and let $\gvct{q}$ be the least right common multiple of $\gvct{q}_1$ and
$\gvct{q}_2$. Then $\Gamma(R_1,R_2)=
\Gamma(R(\gvct{q}_1))\cap\Gamma(R(\gvct{q}_2))=
P(\gvct{q}\gvct{H}+\gvct{q}_1\gvct{r}\gvct{\bar{q}}_2\Z)$.
\end{lemma}
Such an $\gvct{r}$ is by no means unique. Its existence follows from the
uniqueness of the left (or right) prime power decomposition of integer
quaternions. Alternatively we may decompose $\Gamma(R_1,R_2)$ as follows:
\begin{lemma}
Under the conditions of the previous lemma, we have
$\Gamma(R_1,R_2)=P(\gvct{q}\gvct{H}+\gvct{q}_1\gvct{H}\gvct{\bar{q}}_2)=
P(\gvct{q}\gvct{H}+\gvct{q}_2\gvct{H}\gvct{\bar{q}}_1)$.
\end{lemma}
Note that $\gvct{q}\gvct{H}+\gvct{q}_2\gvct{H}\gvct{\bar{q}}_1$ is no ideal
and hence $\Gamma(R_1,R_2)$ is neither an ordinary CSL nor a multiple
of an ordinary CSL. Note further that
$P(\gvct{q}\gvct{H}+\gvct{q}_1\gvct{H}\gvct{\bar{q}}_2)/P(\gvct{q}\gvct{H})$
is a cyclic group of order
$\frac{|\gvct{q}|^2}{\max(|\gvct{q}_1|^2,|\gvct{q}_1|^2)}$ and that
$P(\gvct{q}\gvct{H})$ is a multiple of an ordinary CSL ($\gvct{q}$ is not
primitive here!).
The next lemma tells us under which conditions different pairs of
CSLs give rise to different MCSLs:
\begin{lemma}
Let $\gvct{q}_i$ be primitive quaternions with $|\gvct{q}_i|^2=p^{\alpha_i}$,
where $p$ is a prime and $\alpha_1\geq\alpha_2\geq\alpha_4,
\alpha_3\geq\alpha_4$. Let $\gvct{q}_{ij}$ with $|\gvct{q}_{ij}|^2=p^{\alpha_{ij}}$
be the greatest left common divisor of $\gvct{q}_i$ and $\gvct{q}_j$. Then
$\Gamma(R_1)\cap\Gamma(R_2)=\Gamma(R_3)\cap\Gamma(R_4)$ if and only if
(in case of $\alpha_1=\alpha_2$ possibly after interchanging $\gvct{q}_1$
and $\gvct{q}_2$)
$\alpha_1=\alpha_3,\alpha_2-\alpha_{12}=\alpha_4-\alpha_{34},
\alpha_1-\alpha_{13}\leq\min(\alpha_4-\alpha_{34},\alpha_{34})$ and
$\alpha_2-\alpha_{24}\leq\min(\alpha_4-\alpha_{34},\alpha_{34})$ are satisfied.
\end{lemma}
Thus we can calculate the number $f_2(\Sigma)$ of different MCSLs which are
intersections of at most two ordinary CSLs:
\begin{theo}
Let $p$ be an odd prime number. Then
\BAL
f_2(p^r)=&(r/2+1/2)\,(p+1)p^{r-1}+
(r/2-1)p^{r-2}+
(r/2-[r/2])p^{r-4}\nonumber\\
&\mbox{}+\frac{p^{r-1}-p^{r-2[r/3]-1}}{p^2-1}+
\frac{p^{4[r/3]-r+2}-p^{4[r/2]-r-2}}{2(p^2-1)},
\EAL
where $[x]$ is Gauss' symbol denoting the largest integer $n$ such that
$n\leq x$.
\end{theo} 
Note that $f_2(p^r)=f_m(p^r)$ for $r=1,2$. Thus we know $f_m(\Sigma)$ for all
$\Sigma$ that are free from third powers. The more complex general case will
be presented elsewhere~\cite{pzmcsl2}.

\section*{Acknowledgements}
The author is very grateful to Michael Baake for interesting and stimulating
discussions
on the present subject. Financial support by the Austrian
Academy of Sciences (APART-program) and the
EU Research Training Network
``Quantum Probability with Applications to Physics, Information Theory and
Biology''
is gratefully acknowledged.


\medskip

\begin{thebibliography}{10}

\bibitem{boll70}
Bollmann W 1970
 {\em Crystal Defects and Crystalline Interfaces}
 (Berlin: Springer)

\bibitem{boll82}
Bollmann W
 1982 {\em Crystal Lattices, Interfaces, Matrices}
 published by the author, Geneva.

\bibitem{grim74}
Grimmer H
 1974
 Disorientations and coincidence rotations for cubic lattices
 {\em Acta Cryst.}, A \textbf{30} 685--8

\bibitem{gribo74}
Grimmer H, Bollmann W, and Warrington D~H
 1974
 Coincidence-site lattices and complete pattern-shift lattices in
  cubic crystals
 {\em Acta Cryst.} A \textbf{30} 197--207

\bibitem{ble81}
Bleris G~L and Delavignette P
 1981
 A new formulation for the generation of coincidence site lattices
  ({CSL}'s) in the cubic system
 {\em Acta Cryst.} A \textbf{37} 779--86

\bibitem{grim84}
Grimmer H
 1984
 The generating function for coincidence site lattices in the cubic
  system
 {\em Acta Cryst.} A \textbf{40} 108--12.

\bibitem{baa97}
Baake M
 1997 Solution of the coincidence problem in dimensions $d\leq 4$.
 In {\em {T}he {M}athematics of {L}ong-{R}ange
  Aperiodic {O}rder} ed. R~V Moody (Dordrecht: Kluwer) pp 9--44

\bibitem{pzcsl1}
Zeiner P
 Symmetries of coincidence site lattices of cubic lattices
 {\em Z. Kristallogr.} \textbf{220} 915--25

\bibitem{pzcsl2}
Zeiner P
 Coincidences of hypercubic lattices in 4 dimensions
 {\em Z. Kristallogr.}, in press

\bibitem{plea96}
Pleasants P~A~B, Baake M, and Roth J
 1996
 Planar coincidences for $n$--fold symmetry
 {\em J. Math. Phys.} \textbf{37} 1029--58 (math.MG/0511147)

\bibitem{slo02a}
Diggavi S~N, Sloane N~J~A, and Vaishampayan V~A
 2002
 Asymmetric multiple description lattice vector quantizers
 {\em IEEE Transactions Information Theory} \textbf{48} 174--91

\bibitem{slo02b}
Sloane N~J~A and Beferull-Lozano B
 2003
 Quantizing using lattice intersections
 {\em Discrete and Computational Geometry}
 ed B Aronov, S Basu, J Pach and M Sharir (Berlin: Springer)
 pp 799--824
 (math.CO/0207147)

\bibitem{baagri05}
Baake M and Grimm U
 2005
 Multiple planar coincidences with $N$--fold symmetry
 \textit{Preprint} math.MG/0511306

\bibitem{pzmcsl2}
Zeiner P
 Multiple CSLs for cubic lattices,
 in preparation

\bibitem{koecheng}
Koecher M and Remmert R
 1991 {H}amilton's Quaternions
 In  {\em Numbers} ed. Ebbinghaus H-D et al pp. 155--81 (Springer)

\bibitem{hurw}
Hurwitz A 1919
 {\em {V}or\-le\-sun\-gen \"uber die {Z}ah\-len\-theo\-rie der
  {Q}ua\-ter\-nio\-nen}
 (Berlin: Springer) 

\bibitem{val}
du~Val P 1964
 {\em Homographies, Quaternions and Rotations}
 (Oxford: Clarendon Press) 

\end{thebibliography}


\end{document}